\documentclass[a4paper]{article}
\usepackage{amsmath,amsfonts,amssymb,color,array}
\usepackage{epsfig,cite}
\usepackage{graphicx,epstopdf}
\newtheorem{newexample}{Example}[section]

\newtheorem{newremark}{Remark}[section]

\newtheorem{newkey}{Keywords}
\newenvironment{keywords}
{\begin{newkey}\rm}
	{\end{newkey}}

\newtheorem{newdef}{Definition}[section]

	\newtheorem{theorem}{Theorem}[section]
	\newtheorem{lemma}{Lemma}[section]

\begin{document}

	\title{Exponentially convergent method for time-fractional evolution equation}
	\author{V.Vasylyk\thanks{Institute of
			Mathematics of NAS of Ukraine,
			3 Tereshchenkivs'ka Str., Kyiv-4, 01601, Ukraine
			({\tt vasylyk@imath.kiev.ua}).}
		\and V.L.Makarov\thanks{Institute of
			Mathematics of NAS of Ukraine,
			3 Tereshchenkivs'ka Str., Kyiv-4, 01601, Ukraine
			({\tt makarov@imath.kiev.ua}).}
	}
	
	\date{\null}
	\maketitle

\begin{abstract}
	An exponentially convergent numerical method for solving a differential equation with a right-hand fractional Riemann-Liouville time-derivative and an unbounded operator coefficient in Banach space is proposed and analysed for a homogeneous/inhomogeneous equation of the Hardy-Tichmarsh type. 
	We employ a solution representation by the Danford-Cauchy integral on  hyperbola that envelopes spectrum of the operator coefficient with a subsequent application of an exponentially convergent quadrature. 
	To do that, parameters of the hyperbola are chosen so that the integration function has an analytical extension into a strip around the real axis and then apply the Sinc-quadrature. 
	We show the exponential accuracy and illustrate the results by a numerical example confirming the {\it a priori} estimate. 
	Existence conditions for the solution of the inhomogeneous equation are established. 
\end{abstract}	


\section{Exponentially convergent method for homogeneous Hardy-Titchmarsh--type equation with an unbounded operator coefficient in Banach space}
\subsection{Introduction}
In recent years, interest to differential equations with fractional derivatives has been significantly growing. 
This is due to the fact that the fractional analysis has found out wide application in the modeling of natural and social phenomena. 
The related problems arise most intensively when simulating the linear visco-elastic behavior. 
Fractional derivatives naturally appear in models of anomalous diffusion, control theory, electrodynamics, nonlinear hydroacoustics as well as diffusion which is one of the most significant transport mechanisms in physics. 

The classical diffusion $ \partial_t u-Au = f $ model that adopts the first-order time- derivative, $ \partial_t u $, and the Laplace operator, $ Au = - \Delta u $, is based on assumption that particle motions are Brownian. One of the distinguishing features of the Brownian motions is the linear growth of the mean-square displacement of particles with time $ t $. 

Over the past few decades, a series of experimental studies indicated that assumption on the Brownian type particle motions may not be sufficient for an accurate description of physical processes, and the mean square particle displacement may grow either sublinearly or superlinearly with time $ t $ that is noted in the physical literature as subdiffusion or superdiffusion, respectively \cite{VAS:jlz}. 
These studies cover an rather wide and diverse spectrum of practical applications in engineering, physics, biology, and finance including the electron transport in copier, thermal diffusion in fractal domains, protein transport in cell membranes and so on.

The original equation connects fractional time-derivative and a spatial operator $ A $. 
The input data and the outputs (solutions) are combined using the so-called solving operator $ E_ {1+ \alpha} \left (-At ^ {1+ \alpha} \right) $ which maps the input data to the output. 
A link between the fractional derivative and power of the operator has been examined, for example, in \cite{VAS:ash}. 
In \cite{VAS:g1}, an algorithmic representation of fractional powers of  positive operator $ A $ was proposed.

We consider a differential equation with a Riemann-Liouville fractional derivative on a semi-axis and an unbounded operator coefficient. 
Various methods of discretisation of the differential equations are considered in the literature \cite{VAS:mt}. 
A disadvantage of these discretisations \cite{VAS:mt} is that constants in the accuracy estimates depend exponentially on $ t $. 
Exponential convergence of the discrete approximations to operator-valued functions and differential equations with unbounded operator coefficients plays a crucial role in constructing algorithms of optimal or near-optimal complexity \cite{VAS:GMV-mon,VAS:ghk3}.

\subsection{Presentation of the solution to differential equations with fractional derivatives and unbounded operator coefficients in Banach spaces }

Let us consider the following problem
\begin{equation} \label{VAS:eq-post}
	\begin{split}
		-&{}_{t} D_{\infty }^{\alpha +1} u(t)+Au(t)=0,\quad t\in (0,\infty ),\, \, \, \alpha \in (-1,1), \\[3mm]
		&u(0)=u_{0},
	\end{split}
\end{equation}
where ${}_{t} D_{\infty }^{1+\alpha } $ is the (right) Riemann-Liouville derivative defined by
\begin{equation} \label{VAS:R-L-poh}
	{}_{t} D_{\infty }^{\nu}f(t)=  \begin{cases}\displaystyle
		{\frac{1}{\Gamma (-\nu )} \int _{t}^{\infty }(s-t)^{-\nu -1} f(s)ds ,\qquad \qquad \qquad \nu <0}, \\
		\displaystyle {\frac{1}{\Gamma \left(1-\left\{\nu \right\}\right)} \left(-\frac{d}{dt} \right)^{\left[\nu \right]+1} \int _{t}^{\infty }\frac{f(s)}{(s-t)^{\left\{\nu \right\}} } ds ,\quad \nu \ge 0},
	\end{cases}
\end{equation}
$A$ is a strongly positive operator with everywhere dense  domain $D(A)$ in the Banach space $X$. Its spectrum lies in the sector 
$$ \Sigma (A)=\left\{z=\rho_{0} +re^{i\theta }:\, \, r\in [0,\infty ),\, 0<|\theta |\le \varphi <\frac{1}{2}  {\pi }, \, \rho_{0}>0 \right\}$$
and the following estimate holds true on its boundary $\Gamma _{\Sigma } $ and outside 
\begin{equation}\label{VAS:OcRez}
	\left\| (zI-A)^{-1} \right\| \le {M}/{(1+|z|)} 
\end{equation}
with a positive constant $M$.

The following theorem was proven in \cite{VAS:gmv-HT}.
\begin{theorem} Solution to equation \eqref{VAS:eq-post} is 
	\begin{equation} \label{VAS:i2-3}
		u(t)=\exp (-A^{1/(1+\alpha )} t)u(0),
	\end{equation}
	when the following limits 
	\begin{equation} \label{VAS:i2-2}
		{\mathop{\lim }\limits_{s\to \infty }} \left[(s-t)^{\alpha +1} {}_{s} D_{\infty }^{\alpha } u(s)\right]=0, \
		{\mathop{\lim }\limits_{s\to \infty }} \left[(s-t)^{\alpha } {}_{s} D_{\infty }^{\alpha -1} u(s)\right]=0
	\end{equation}
	are fulfilled.   
	It means that the solution operator takes the form
	\begin{equation} \label{VAS:i2-3-1}
		S(t,A)=\exp (-A^{1/(1+\alpha )} t).
	\end{equation}
\end{theorem}

The formula \eqref{VAS:i2-3} was derived by using operator (\ref{VAS:R-L-poh}) in equation (\ref{VAS:eq-post}) and accounting for (\ref{VAS:i2-2}) so that we arrive at 
\begin{equation} \label{VAS:i2-3-2}
	-u(t)+A\, {}_{t} D_{\infty }^{-(\alpha +1)} u(t)=0,\quad t\in (0,\infty ),\, \, \, \alpha \in (-1,1),
\end{equation}
that coincides with the Hardy-Titchmarsh integral equation \cite{VAS:ht} with  an operator coefficient (when changing $ \alpha $ to $ \alpha + 1 $). 
Therefore, its solution possesses  \eqref{VAS:i2-3}.

Let us construct a solution to equation  \eqref{VAS:eq-post} using the Dunford-Cauchy integral. 
As it was shown in \cite{VAS:gmv-HT} (see also \cite{VAS:GMV-mon}), the solution can be written as 
\begin{equation}\label{VAS:Zobr}
	u(t) = \frac{1}{2\pi i} \int\limits_\Gamma  {\rm e}^{ - t{\kern 1pt} z^{1 /(1 + \alpha )}} 
	\left[ {{{\left( {zI - A} \right)}^{ - 1}} - \frac{1}{z}I} \right]u(0)dz ,
\end{equation}
where $\Gamma$ is a smooth curve that envelopes spectrum of operator $A$.

Let us call
\begin{multline}\label{VAS:spectr-hyp}
	\Gamma_0=\left\{ z\left( s\right) =a_0\left( \cosh s -1\right) +\rho_{0}-ib_0\sinh s: \, s\in \left( -\infty,\infty\right),  \right. \\
	\left. b_0=a_0\tan \varphi, \,\, a_0: \, \rho_{0}-a_0>0 \right\}
\end{multline}
the spectral hyperbola which envelopes from the left the spectrum $\Sigma(A)$. 
The hyperbola asymptotes for  $s\to \pm\infty$ form angle $\varphi$ with the real axis. 
To construct an exponentially convergent numerical method, it is necessary to choose the contour of integration so that it envelopes $\Gamma_0$, bypassing it from the left and at the same time an integral \eqref{VAS:Zobr} exists \cite{VAS:GMV-mon}. We will look for such a contour of hyperbolic shape
\begin{equation}\label{VAS:int-hyp}
	\Gamma_I=\left\{ z\left( s\right) =a_I\left( \cosh s -1\right) +q - ib_I\sinh s: \, s\in \left( -\infty,\infty\right) \right\}
\end{equation}
and call it the integral hyperbola.
Having chosen $\Gamma_I$ as integration contour, provided by \eqref{VAS:Zobr}, the solution representation takes the form
\begin{equation}\label{VAS:param-zobr}
	u(t) = \frac{1}{2\pi i} \int\limits_{ - \infty }^\infty  {{F_A}(t,\xi )u(0) d\xi },
\end{equation}
where
\[{F_A}(t,\xi ) = {{\rm e} ^{ - z{{(\xi )}^{1/(1 + \alpha )}}t}}({a_I}\sinh\xi  - i\,\,{b_I} \cosh\xi )\left[ {{{\left( {z(\xi )I - A} \right)}^{ - 1}} - \frac{1}{{z(\xi )}}I} \right],\]
\[z(\xi ) = a_I\left( \cosh \xi -1\right) +q - ib_I\sinh \xi .\]

\subsection{Study of parameters of integration contours} 

By Theorem 2.2 from \cite{VAS:GMV-mon} coefficients of the integral hyperbola should be chosen in such way that:\\
1.~Operator-valued function  ${F_A}(t,\xi )\in H^p(D_d),$ $\forall t\ge0$ (see \cite{VAS:stenger}) so that 
\begin{equation*}
	D_d=\{z \in \mathbb{C}: - \infty < \Re z < \infty, |\Im z|<d \},
\end{equation*}
\begin{equation*}
	\Vert{\mathcal{F}}\Vert_{\mathbf{H}^{p}(D_{d})}:=
	\begin{cases}
		\lim_{\epsilon\rightarrow 0}\left(  \int_{\partial D_{d}(\epsilon)}
		\Vert\mathcal{F}(z)\Vert^{p}|dz|\right)  ^{1/p} & \text{if $1\leq p<\infty$},\\
		\lim_{\epsilon\rightarrow 0}\sup_{z\in D_{d}(\epsilon)}\Vert\mathcal{F}(z)\Vert & \text{if $p=\infty$}.
	\end{cases}
\end{equation*}
$D_d(\epsilon)$ is defined for $0< \epsilon <1$ as follows
$$ 
D_d(\epsilon)=\{z \in \mathbb{C}: | \Re z| < 1/\epsilon, |\Im
z|<d(1-\epsilon)\}.
$$ 
2.~There exist positive constants  $c,$ $\delta,$ such that the following estimate 
$$\left\| {F_A}(t,\xi )\right\| \le c {\rm e}^{-\delta \left| \xi\right| }, \quad \xi\in(-\infty,\infty), \, t\ge 0$$
holds true.

\medskip 

Let us focus on the first condition by constructing a parametric family of curves
\[
\Gamma_I(\nu)=z(s+i\nu), \quad \nu\in \left[ -\frac{d_1}{2},\frac{d_1}{2}\right].
\]
Then
\[
z(s+i\nu)= a(\nu) \cosh s -a_I+q - ib(\nu)\sinh s,
\]
with
\[
\begin{split}
	& a(\nu)=a_I \cos \nu+b_I \sin \nu=\sqrt{a_I^2+b_I^2}\sin(\psi+\nu), \\
	& b(\nu)=b_I \cos \nu-a_I \sin \nu=\sqrt{a_I^2+b_I^2}\cos(\psi+\nu),\\
	&  \cos(\psi)=\frac{b_I}{\sqrt{a_I^2+b_I^2}}.
\end{split}
\]

Let us choose coefficients $a_I, \, b_I$ in such a way that curve  $\Gamma_I(\nu)$  becomes a straight line parallel to the imaginary axis for $\nu= -\frac{d_1}{2}$ and coincided with the spectral hyperbola as $\nu= \frac{d_1}{2}$. The requirements lead to obvious equalities
\[
\displaystyle a\left( -\frac{d_1}{2}\right) =0;\ \ 
\displaystyle a\left( \frac{d_1}{2}\right) =a_0;\ \ 
\displaystyle b\left( \frac{d_1}{2}\right) =b_0,
\]
i.e.,
\[
\begin{cases}
	& \displaystyle a_I\cos \frac{d_1}{2}-b_I \sin \frac{d_1}{2} =0,\\
	& \displaystyle -a_I\sin \frac{d_1}{2}+b_I \cos \frac{d_1}{2} =b_0,\\
	& \displaystyle a_I\cos \frac{d_1}{2}+b_I \sin \frac{d_1}{2} =a_0.
\end{cases}
\]
and, therefore,
\[
\displaystyle a_I=a_0\cos \frac{d_1}{2}-b_0 \sin \frac{d_1}{2}; \ \
\displaystyle b_I=a_0\sin \frac{d_1}{2}+b_0 \cos \frac{d_1}{2};\ \
\displaystyle 2a_I\cos \frac{d_1}{2} =a_0
\]
whose first and third equations deduce  
\[
\cos d_1=\tan \varphi \sin d_1.
\]
As a consequence, 
\begin{equation*}
	\begin{gathered}
		\displaystyle d_1=\frac{\pi}{2}-\varphi;\ \
		\displaystyle a_I=a_0 \frac{\cos\left( \frac{\pi}{4}+\frac{\varphi}{2}\right)} {\cos\varphi}; \ \
		\displaystyle b_I=a_0 \frac{\sin\left( \frac{\pi}{4}+\frac{\varphi}{2}\right) }{\cos\varphi},\\
		\displaystyle a(\nu)=a_0 \frac{\cos\left( \frac{\pi}{4}+\frac{\varphi}{2}-\nu \right)} {\cos\varphi};\ \ 
		\displaystyle b(\nu)=a_0 \frac{\sin\left( \frac{\pi}{4}+\frac{\varphi}{2}-\nu \right) }{\cos\varphi}.
	\end{gathered}
\end{equation*}

Next, we need to define the parameter $q$. For that, we assume that  
\begin{enumerate}
	\item Curve $\displaystyle \Gamma(\nu), \, \forall \nu\in \left[ -\frac{d_1}{2},\frac{d_1}{2}\right]$ is in the right half-plane and  ${F_A}(t,\xi )\in H^p(D_d)$, $\forall t\ge0$.
	\item Curve $\Gamma(\nu)$ does not intersect spectrum of the operator $A$.
\end{enumerate} 
These two conditions yield 
\[
q-a_I>0;\ \
q-a_I<\rho_{0}
\]
which lead to  
$
a_I<q<\rho_{0}+a_I
$
or, remembering previous results, to 
\[ \displaystyle
a_0 \frac{\cos\left( \frac{\pi}{4}+\frac{\varphi}{2}\right)} {\cos\varphi}<q< a_0 \frac{\cos\left( \frac{\pi}{4}+\frac{\varphi}{2}\right)} {\cos\varphi}+\rho_0.
\]

Let us choose 
\[ \displaystyle
a_0=\frac{\cos\varphi}{\cos\left( \frac{\pi}{4}+\frac{\varphi}{2}\right)}; \ \ \ 
q=\frac{\rho_{0}}{2}+1.
\]
Then, obviously, the specified conditions will be fulfilled. Thus, we proved the following lemma.

\begin{lemma}\label{VAS:Lem1} 
	Let parameters of curves $\Gamma_0$ and $\Gamma_I$ be
	\begin{equation}\label{VAS:param_hyp}
		\begin{split}
			& \displaystyle a_0=\frac{\cos\varphi}{\cos\left( \frac{\pi}{4}+\frac{\varphi}{2}\right)}; \quad b_0=\frac{\sin\varphi}{\cos\left( \frac{\pi}{4}+\frac{\varphi}{2}\right)}, \\
			& \displaystyle a_I=1; \quad b_I=\tan\left( \frac{\pi}{4}+\frac{\varphi}{2}\right); \quad
			d_1=\frac{\pi}{2}-\varphi; \quad q=\frac{\rho_{0}}{2}+1.
		\end{split}
	\end{equation}
	Then the family $\Gamma_I(\nu)$ for $\displaystyle \nu\in \left[ -\frac{d_1}{2},\frac{d_1}{2}\right]$ will be in the right half-plane and will not cross the spectral hyperbola $\Gamma_{\Sigma }$.
\end{lemma}

Accounting for 
\[
\sqrt{a_I^2+b_I^2}=\sqrt{1+\tan^2\left( \frac{\pi}{4}+\frac{\varphi}{2}\right)} 
=\frac{1}{ \left|\cos \left( \frac{\pi}{4}+\frac{\varphi}{2} \right)\right| },
\]
we get
\[
\sin(\psi)=\frac{1}{\sqrt{a_I^2+b_I^2}}=\cos \left( \frac{\pi}{4}+\frac{\varphi}{2} \right)=\sin \left( \frac{\pi}{4}-\frac{\varphi}{2} \right);\ \ 
\cos ( \psi )=\sin \left( \frac{\pi}{4}+\frac{\varphi}{2} \right).
\]

Next, we examine when ${F_A}(t,\xi )\in H^p(D_d),$ $\forall t\ge0$ and exponentially decays as $\xi \to \pm\infty$. For $z^{\frac{1}{1+\alpha}}(\xi)$, 
\[
\left[ z\left(\xi+i\nu \right) \right]^{\frac{1}{1+\alpha}}= \left[ r\left(\xi,\nu \right) \right]^{\frac{1}{1+\alpha}} \exp\left\lbrace i\frac{\theta(\psi+\nu,\xi)}{1+\alpha} \right\rbrace ,
\]
where
\[
r\left(\xi,\nu \right)=\left| z(s+i\nu)\right| , 
\]
\[
\begin{split}
	&\cos \theta(\psi+\nu,\xi)=\frac{a(\nu)\cosh\xi+\frac{\rho_0}{2}}{r(\xi,\nu)}= \frac{\sqrt{a_I^2+b_I^2}\sin(\psi+\nu)\cosh\xi+\frac{\rho_0}{2}}{r(\xi,\nu)},\\
	&\sin \theta(\psi+\nu,\xi)= \frac{-\sqrt{a_I^2+b_I^2}\cos(\psi+\nu)\sinh\xi}{r(\xi,\nu)}.
\end{split}
\]

Then, to fulfil
\[
\left| \exp\left\lbrace -z^{\frac{1}{1+\alpha}}(\xi) t \right\rbrace \right| \to 0, \quad \xi\to\pm\infty,
\]
it would be
\[
\cos \frac{\theta(\psi+\nu,\xi)}{1+\alpha} >0.
\]
Accounting for that  $1+\alpha>0$, one concludes 
\begin{equation}\label{VAS:for-theta}
	\left| \frac{\theta(\psi+\nu,\xi)}{1+\alpha} \right| <\frac{\pi}{2};\ \ \ 
	\left| \theta(\psi+\nu,\xi) \right| <\frac{\pi}{2}(1+\alpha).
\end{equation}

Let us consider the case $-1<\alpha\le 0$. Then \eqref{VAS:for-theta} deduces 
\[
\left| \sin \theta(\psi+\nu,\xi) \right| <\sin\frac{\pi}{2}(1+\alpha)=\cos\frac{\pi\alpha}{2}.
\]
Therefore, we arrive at the following condition
\[
\left| \frac{-\sqrt{a_I^2+b_I^2}\cos(\psi+\nu)\sinh\xi}{r(\xi,\nu)}\right| < \cos\frac{\pi\alpha}{2}
\]
where, because 
\[
\psi= \frac{\pi}{4}-\frac{\varphi}{2};\quad \nu\in \left[ -\frac{d_1}{2}, \frac{d_1}{2}\right]= \left[ -\frac{\pi}{4}+\frac{\varphi}{2}, \frac{\pi}{4}-\frac{\varphi}{2}\right], 
\]
evaluating the left hand side of the inequality deduces 
\[
\psi+\nu \in \left[ 0, \frac{\pi}{2}- \varphi \right]
\]
and, as a result, 
\begin{multline*}
	\frac{\sqrt{a_I^2+b_I^2}\left|\sinh\xi\right|}{r(\xi,\nu)}\\
	= \frac{\sqrt{a_I^2+b_I^2}\left|\sinh\xi\right|}{\left[ \left( \sqrt{a_I^2+b_I^2}\sin(\psi+\nu)\cosh\xi+ \frac{\rho_0}{2} \right)^2 + \left( \sqrt{a_I^2+b_I^2}\cos(\psi+\nu)\sinh\xi \right)^2 \right] ^{0.5}} \\
	=\frac{\left|\tanh\xi\right|}{\left[ \left( \sin(\psi+\nu) + \frac{\rho_0}{2\sqrt{a_I^2+b_I^2} \cosh\xi} \right)^2 + \left( \cos(\psi+\nu)\tanh\xi \right)^2 \right] ^{0.5}} \\
	\le\frac{\left|\tanh\xi\right|}{\left[ \sin^2(\psi+\nu) + \cos^2(\psi+\nu)\tanh^2\xi \right] ^{0.5}}.
\end{multline*}

Function $f(y)=y/\sqrt{a^2+b^2y^2},\, y\in [0,\infty)$ is monotonically increasing alike $\tanh( y)$, thus, 
\[
\frac{\left|\tanh\xi\right|}{\left[ \sin^2(\psi+\nu) + \cos^2(\psi+\nu)\tanh^2\xi \right] ^{0.5}}\le \frac{1}{\left[ \sin^2(\psi+\nu) + \cos^2(\psi+\nu) \right] ^{0.5}}=1
\]
and 
\[
\frac{\sqrt{a_I^2+b_I^2}\left|\cos(\psi+\nu)\sinh\xi\right|}{r(\xi,\nu)} \le \left| \cos(\psi+\nu) \right|.
\]

This means that if 
\begin{equation}\label{VAS:obm}
	\left| \cos(\psi+\nu) \right| < \cos\frac{\pi\alpha}{2},
\end{equation}
then the inequalities \eqref{VAS:for-theta} hold true and, therefore, when $-1<\alpha\le 0$,
\[
\psi+\nu>-\frac{\pi\alpha}{2},
\]
or
\begin{equation}\label{VAS:um-phi}
	\varphi<(1+\alpha)\frac{\pi}{2}.
\end{equation}

In the case $0\le\alpha<1$, 
\[
\cos\theta(\psi+\nu,\xi)=\frac{\cos^{-1} \left( \frac{\pi}{4}+\frac{\varphi}{2} \right)\sin(\psi+\nu) \cosh\xi+\frac{\rho_0}{2}}{r(\xi,\nu)}>0,
\]
provided by 
\[
\psi+\nu\in\left[ 0, \frac{\pi}{2}- \varphi  \right] \Rightarrow \sin(\psi+\nu)\ge0,
\]
\[
0< \frac{\pi}{4}+\frac{\varphi}{2}< \frac{\pi}{2} \Rightarrow \cos\left( \frac{\pi}{4}+\frac{\varphi}{2} \right)>0
\]
and, as a result, 
\[
\left| \theta(\psi+\nu,\xi) \right| <\frac{\pi}{2}<\frac{\pi}{2}(1+\alpha).
\]

Let us find $\left\| F_A(t,w)u(0)\right\| $ in the strip $w=\xi+i\nu \in D_1.$ According to Theorem 2.4 from \cite{VAS:GMV-mon} for $m=0$,
\begin{equation}\label{VAS:oc11}
	\left\| \left[ \left( z(\xi)I-A\right)^{-1}-\frac{I}{z(\xi)} \right] u_0\right\| \le \frac{1}{\left| z(w) \right| }\frac{(1+M)K}{\left( 1+ \left| z(w) \right|\right) ^\gamma} \left\| A^\gamma u(0) \right\|, \ u(0)\in D(A^\gamma) ,
\end{equation}
where the positive constant $K$ depends on  $\alpha$.

Furthermore,
\begin{multline*}
	\left| \frac{z^\prime(w)}{z(w)}\right| =\left| \frac{a(\nu)\sinh\xi-ib(\nu)\cosh\xi}{a(\nu)\cosh\xi +\frac{\rho_0}{2}-ib(\nu)\sinh\xi} \right| \\
	=\left[ \frac{a^2(\nu)\sinh^2\xi+ b^2(\nu)\cosh^2\xi}{\left( a(\nu) \cosh\xi +\frac{\rho_0}{2}\right)^2  +b^2(\nu)\sinh^2\xi} \right]^{0.5}\\
	= \left[ \frac{a^2(\nu)\tanh^2\xi+ b^2(\nu)}{\left( a(\nu) +\frac{\rho_0}{2\cosh\xi}\right)^2  +b^2(\nu)\tanh^2\xi} \right]^{0.5} \\
	\le \left[ \frac{a^2(\nu)\tanh^2\xi+ b^2(\nu)}{ a^2(\nu) + b^2(\nu)\tanh^2\xi} \right]^{0.5}\le \max \left\lbrace 1,\left| \frac{a(\nu)}{b(\nu)} \right| \right\rbrace = \max \left\lbrace 1,\tan(\psi+\nu) \right\rbrace.
\end{multline*}

Accounting for $\displaystyle \psi+\nu\in \left[ 0, \frac{\pi}{2}- \varphi  \right] \Rightarrow \tan(\psi+\nu)<\infty,$ we get
\begin{equation}\label{VAS:oc-z'}
	\left| \frac{z^\prime(w)}{z(w)}\right|\le \max \left\lbrace 1,\tan\left( \frac{\pi}{2}- \varphi\right)  \right\rbrace=c_1<\infty
\end{equation}
and
\begin{multline*}
	\left| \exp\left\lbrace -z^{\frac{1}{1+\alpha}}(\xi+i\nu) t \right\rbrace  \right| = \left| \exp\left\lbrace  -\left[ r\left(\xi,\nu \right) \right]^{\frac{1}{1+\alpha}} \exp\left\lbrace i\frac{\theta(\psi+\nu,\xi)}{1+\alpha} \right\rbrace t \right\rbrace\right| \\
	= \exp\left\lbrace  - \left[ r\left(\xi,\nu \right) \right]^{\frac{1}{1+\alpha}} \cos \frac{\theta(\psi+\nu,\xi)}{1+\alpha} t \right\rbrace \le 1,
\end{multline*}
when condition \eqref{VAS:um-phi} is true for $-1<\alpha\le0$ and $\varphi<\frac{\pi}{2}$ when $0<\alpha<1.$

If $\nu\in \left( -\frac{d_1}{2},\frac{d_1}{2}\right] ,$ then
\begin{multline*}
	\frac{1}{\left( 1+\left| z(w) \right| \right)^\gamma }\le \frac{1}{\left[ \left( \frac{\sin(\psi+\nu)}{\cos\left( \frac{\pi}{4}+\frac{\varphi}{2}\right)}\cosh\xi+ \frac{\rho_0}{2} \right)^2 + \left( \frac{\cos(\psi+\nu)}{\cos\left( \frac{\pi}{4}+\frac{\varphi}{2}\right)}\sinh\xi \right)^2 \right] ^{0.5\gamma}} \\
	= \frac{1}{\cosh^\gamma\xi \left[ \left( \frac{\sin(\psi+\nu)}{\cos\left( \frac{\pi}{4}+\frac{\varphi}{2}\right)}+ \frac{\rho_0}{2\cosh\xi} \right)^2 + \left( \frac{\cos(\psi+\nu)}{\cos\left( \frac{\pi}{4}+\frac{\varphi}{2}\right)}\tanh\xi \right)^2 \right] ^{0.5\gamma}}\\
	\le \frac{1}{C_1 \cosh^\gamma\xi }\le C_2{\rm e}^{-\gamma|\xi|},
\end{multline*}
where 
\[
C_1=\min_{\xi,\nu} \left[ \left( \frac{\sin(\psi+\nu)}{\cos\left( \frac{\pi}{4}+\frac{\varphi}{2}\right)}+ \frac{\rho_0}{2\cosh\xi} \right)^2 + \left( \frac{\cos(\psi+\nu)}{\cos\left( \frac{\pi}{4}+\frac{\varphi}{2}\right)}\tanh\xi \right)^2 \right] ^{0.5\gamma}>0.  
\]

For $\nu=-d_1/2$,
$
a(\nu)=0, b(\nu) =\cos^{-1}\left( \frac{\pi}{4}+\frac{\varphi}{2}\right)
$
and, therefore, 
\begin{multline*}
	\frac{1}{\left( 1+\left| z(w) \right| \right)^\gamma }= \frac{1}{\left( 1+ 
		\left[ \displaystyle \frac{\rho_0^2}{4} +\frac{\sinh^2\xi}{\cos^2 \left( \frac{\pi}{4}+ \frac{\varphi}{2}\right)} \right]^{0.5} \right) ^\gamma} \le 
	\frac{1}{\left( 1+ \displaystyle \frac{\left| \sinh\xi\right| }{\cos \left( \frac{\pi}{4}+ \frac{\varphi}{2}\right)} \right) ^\gamma} \\
	= \frac{1}{\cosh^\gamma\xi \left( \displaystyle \frac{1}{\cosh\xi}+ \frac{\left| \tanh\xi\right| }{\cos \left( \frac{\pi}{4}+ \frac{\varphi}{2}\right)} \right) ^\gamma}\le \frac{\cos^\gamma \left( \frac{\pi}{4}+ \frac{\varphi}{2}\right)}{\cosh^\gamma\xi} \le \cos^\gamma \left( \frac{\pi}{4}+ \frac{\varphi}{2}\right) {\rm e}^{-\gamma|\xi|}. 
\end{multline*}
Here we have used that function $f(\xi):$ 
\[
f(\xi)= \frac{1}{\cosh \xi}+ \frac{\left| \tanh \xi\right| }{\cos \left( \frac{\pi}{4}+ \frac{\varphi}{2}\right)},
\]
reaches its minimum at $\xi=\pm\infty$, thus, 
\begin{equation}\label{VAS:oc13}
	\frac{1}{\left( 1+\left| z(w) \right| \right)^\gamma }\le C {\rm e}^{-\gamma|\xi|}, \qquad \forall \, \nu \in \left[ -\frac{d_1}{2},\frac{d_1}{2}\right] .
\end{equation}

Accounting for \eqref{VAS:oc11}--\eqref{VAS:oc13}, one gets 
\begin{equation}\label{VAS:ocFa-2}
	\begin{split}
		\left\| {F_A}(t,w ) \right\| \le (1+M)KC\tan\left( \frac{\pi}{2}-\varphi \right)\left\| A^\gamma u(0)\right\| {\rm e}^{-\gamma|\xi|}, \\
		w \in D_1=\left\lbrace \xi+i\nu: \, \xi\in (-\infty,\infty),\,\, \nu \in \left[ -\frac{d_1}{2},\frac{d_1}{2}\right] \right\rbrace .
	\end{split}
\end{equation}	

In summary, we proved the following theorem.

\begin{theorem}\label{VAS:th1}
	Let $A$ be a strongly positive operator with an everywhere dense domain $D(A)$ in the Banach space $X$, with a spectrum in the sector $\Sigma (A)$ and the resolvent  estimate \eqref{VAS:OcRez}. If $\alpha\in(- 1,0)$, then the spectral angle $\varphi$ additionally satisfies the constraint 
	\begin{equation}\label{VAS:um-phi-2}
		\varphi< \frac{ 1}{2} (1+\alpha) {\pi}.
	\end{equation}
	In addition, let the conditions of lemma \ref{VAS:Lem1} and $u(0)\in D\left( A^\gamma\right) $ be fulfilled. 
	
	Then function $ {F_A}(t,\xi )$ has an analytical continuation in the strip $D_1$ and the following estimate 
	\begin{equation}\label{VAS:ocFa}
		\begin{split}
			\left\| {F_A}(t,w ) \right\| \le (1+M)KC\tan\left( \frac{ 1}{2} {\pi}-\varphi \right)\left\| A^\gamma u(0)\right\| {\rm e}^{-\gamma|\xi|}, \\
			w \in D_1=\left\lbrace \xi+i\nu: \, \xi\in (-\infty,\infty),\,\, \nu \in \left[ -\frac{1}{2}  {d_1},\frac{1}{2}  {d_1}\right] \right\rbrace
		\end{split}
	\end{equation}	
	holds true.
\end{theorem}

From the Theorem \ref{VAS:th1} we have that for $ {F_A}(t,\xi )$ the following estimate 
$$ 
\Vert{{F_A}(t,w )}\Vert_{\mathbf{H}^{1}(D_{d})} \le \frac{2C(\varphi,\gamma)}{\gamma} \left\| A^\gamma u(0)\right\|,
$$
$$
C(\varphi,\gamma)= (1+M)KC\tan\left( \frac{1}{2}  {\pi}-\varphi \right)
$$
is fulfilled.

\subsection{Numerical method for the homogeneous equation}\label{VAS:metod-odn} 

Let us use Sinc-quadrature (see  \cite{VAS:GMV-mon,VAS:stenger}) to approximate integral \eqref{VAS:param-zobr}:
\begin{equation}\label{VAS:h-nab}
	u(t) \approx u_{h,N}(t)=\frac{h}{2 \pi i}\sum_{k=-N}^{N} F_A(t,kh)
\end{equation}
with an error $\eta_N(F_A,h)$ such that
$$\|\eta_N(F_A,h)\|=\|u(t)-u_{h,N}(t)\| $$
\[
\le \left\|u_h(t)-\frac{h}{2 \pi i}\sum_{k=-\infty}^{\infty}F_A(t,kh) \right\| + \left\|\frac{h}{2 \pi i}\sum_{|k|>N}F_A(t,kh) \right\|.
\]
For the first term we use the estimate from \cite{VAS:GMV-mon}
\[
\left\|u(t)-\frac{h}{2 \pi i}\sum_{k=-\infty}^{\infty}F_A(t,kh) \right\| \le \frac{1}{2 \pi}\frac{{\rm e}^{-\pi d_1/h}}{2 \sinh{(\pi d_1/h)}} \left\| F_A(t,w) \right\|_{{\bf H}^1(D_{d_1})}.
\]
For the second term from  \eqref{VAS:ocFa} we get
\begin{multline*}
	\left\|\frac{h}{2 \pi i}\sum_{|k|>N}F_A(t,kh) \right\| \le \frac{C(\varphi,\gamma)h \|A^{\gamma}u(0)\| }{2\pi} \sum_{k=N+1}^{\infty}{\rm e}^{-\gamma kh}\\
	\le \frac{C(\varphi,\gamma)h \|A^{\gamma}u(0)\| }{2\pi} {\rm e}^{-\gamma N h}\sum_{k=1}^{\infty} {\rm e}^{-\gamma kh}= \frac{C(\varphi,\gamma)h \|A^{\gamma}u(0)\| }{2\pi}  \frac{{\rm e}^{-\gamma N h}}{{\rm e}^{\gamma h}- 1 }\\
	=\frac{C(\varphi,\gamma)h \|A^{\gamma}u(0)\| }{2\pi}  \frac{{\rm e}^{-\gamma N h}}{\gamma h{\rm e}^{\theta \gamma h} }= \frac{C(\varphi,\gamma) \|A^{\gamma}u(0)\| }{2\pi\gamma {\rm e}^{\theta \gamma h}}  {\rm e}^{-\gamma N h}, \qquad \theta\in (0,1).
\end{multline*}
Here we used Lagrange's theorem on the increment of a function. 
The final estimate takes the form
$$
\|\eta_N(F_A,h)\| \le C_1(\varphi,\gamma) \|A^{\gamma}u(0)\| \left\{{{\rm e}^{-\pi d_1/h}}{\sinh^{-1}{(\pi d_1/h)}} + {\rm e}^{-\gamma N h} \right\},
$$
where $ C_1(\varphi,\gamma)>0$ does not depend on $h,N,t.$ 
Equalising both exponents gives 
$
{2 \pi d_1}/{h}=\gamma Nh,
$
or after transformation, 
\begin{equation}\label{VAS:st-size}
	h=\sqrt{{2 \pi d_1}/{(\gamma N)}}.
\end{equation}

Choosing the step size $h$ in this way, the error of the quadrature formula will satisfy the estimate 
\begin{equation}\label{VAS:na24}
	\|\eta_N(F_A,h)\|\le C_1 \text{exp}{\left(-\sqrt{\frac{1}{2}  {\pi d_1 \gamma} N}\right)} \|A^{\gamma}u(0)\|,
\end{equation}
where constant $c>0$ does not depend on $t,N.$ 

Hence. we proved the following theorem:

\begin{theorem}\label{VAS:Th-OcAlg}
	Let the conditions of Theorem \ref{VAS:th1} be fulfilled. Then for approximate solution of \eqref{VAS:eq-post} by  \eqref{VAS:h-nab}, the accuracy estimate \eqref{VAS:na24} is true if the step $h$ is defined by \eqref{VAS:st-size}.
\end{theorem}

\subsection{Numerical example} 

Consider \eqref{VAS:eq-post} with $\alpha  = \frac{1}{2}$ and 
$$A =  - \frac{\partial ^2}{\partial {x^2}}; \ \ 
D\left( A \right) = \left\{ v\left( x \right) \in \text{H}^1\left( 0,1 \right):\; v\left( 0 \right) = v\left( 1 \right) = 0 \right\},$$
where ${\text{H}}^1\left( {0,1} \right)$ is the Sobolev space, and the initial condition is 
${u_0} = \sin (\pi {\kern 1pt} x).$

The problem has the exact solution 
$u\left( {x,t} \right) = {\operatorname{e} ^{ - {\pi ^4}t}}{\kern 1pt} \sin (\pi {\kern 1pt} x).$
The approximate solution \eqref{VAS:h-nab} of the problem \eqref{VAS:eq-post} is characterised by the error at the point  $t = \pi^{-2},$ $x = 1/2$ which is presented in Table~\ref{VAS:tab1}. 
It shows that the error decreases according to the {\it a priori} estimate \eqref{VAS:na24}.

\begin{table}[!tb]
	\caption{The approximation error of the numerical solution at a point $t = \frac{1}{\pi ^2},$ $x = \frac{1}{2}$.}\label{VAS:tab1}
	\begin{center}
		\begin{tabular}{r@{\quad}l}
			\hline
			$N$    & $\varepsilon_N$ \\
			\hline
			$4$    &   $0.00225$ \\
			$8$    &   $0.00017$ \\
			$16$   &   $7.48280 \cdot 10^{-6} $ \\
			$32$   &   $3.72857 \cdot 10^{-7}$ \\
			$64$   &   $2.29011 \cdot 10^{-9} $ \\
			$128$  &   $1.28593 \cdot 10^{-12} $ \\
			$256$  &   $1.03171 \cdot 10^{-17} $ \\
			$512$  &   $2.94399 \cdot 10^{-23} $ \\
			$1024$ &   $1.25671 \cdot 10^{-30 } $ \\
			\hline
		\end{tabular}
	\end{center}
\end{table}

\section{Exponentially convergent method for the  inhomogeneous Hardy--Titchmarsh type equation with unbounded operator coefficient in the Banach space}

Consider inhomogeneous fractional equation with unbounded operator coefficient in the Banach space of the Hardy--Titchmarsh type, 
\begin{equation}\label{VAS:S-B1.1}
	- {}_tD_\infty ^{\alpha  + 1}u(t) + Au(t) = f\left( t \right),\ t \in (0,\infty ),\,\alpha  \in ( - 1,1); \
	u(0) = {u_0}, 
\end{equation}
where 
${}_tD_\infty ^{1 + \alpha }$ is the right Riemann-Liouville fractional derivative defined in \eqref{VAS:R-L-poh}.
$A$ is a strongly positive operator with everywhere dense  domain $D(A)$ in the Banach space $X$ whose spectrum lies in the sector
$$ \Sigma (A)=\left\{z=\rho_{0} +re^{i\theta }:\, \, r\in [0,\infty ),\, 0<|\theta |\le \varphi <\frac{1}{2}  {\pi }, \, \rho_{0}>0 \right\}$$
and the estimate \eqref{VAS:OcRez} holds true on its boundary $\Gamma _{\Sigma } $ and outside 
\[\left\| {{{(zI - A)}^{ - 1}}} \right\| \leqslant \frac{M}{{1 + |z|}},\] 
with some positive constant $M.$ 

Remark that definition \eqref{VAS:R-L-poh} leads to 
\[
\displaystyle
{}_tD_\infty ^{\alpha  + 1}f(t) = \left\{ {\begin{array}{*{20}{c}}
		{ \displaystyle \frac{{ - 1}}{{\Gamma ( - \alpha )}}\frac{d}{{dt}}\int\limits_t^\infty  {{{\left( {s - t} \right)}^{ - \alpha  - 1}}f(s)ds} ,\quad  - 1 < \alpha  < 0}, \\ 
		{ \displaystyle \frac{1}{{\Gamma \left( {1 - \alpha } \right)}}\frac{{{d^2}}}{{d{t^2}}}\int\limits_t^\infty  {{{\left( {s - t} \right)}^{ - \alpha }}f(s)ds} ,\quad 0 \leqslant \alpha  < 1} ,
\end{array}} \right.
\] 
because $\left[ \alpha  + 1 \right] = 0,$ $\left\{ \alpha  + 1 \right\} = \alpha  + 1,$ for $ - 1 < \alpha  < 0,$ and  
$\left[ {\alpha  + 1} \right] = 1,$ $\left\{ {\alpha  + 1} \right\} = \alpha,$ for $ 0 \leqslant \alpha  < 1.$

\subsection{Solution representation of the inhomogeneous Hardy-Titchmarsh--type equation} 

Let us transform the problem  \eqref{VAS:S-B1.1} to an integral equation. First, we consider the case $ - 1 < \alpha  < 0$ and apply  ${}_tD_\infty ^{ - 2 - \alpha }$ to the left and right hand-sides of \eqref{VAS:S-B1.1}. For the first term, we will use the transformation
\begin{multline*}
	{}_tD_\infty ^{ - 2 - \alpha }{}_tD_\infty ^{1 + \alpha }u(t) \\
	=\frac{{ - 1}}{{\Gamma \left( {2 + \alpha } \right)\Gamma \left( { - \alpha } \right)}}\int\limits_t^\infty  {{{\left( {s - t} \right)}^{1 + \alpha }}\frac{d}{{ds}}\left( {\int\limits_s^\infty  {{{\left( {{s_1} - s} \right)}^{ - 1 - \alpha }}u\left( {{s_1}} \right)d{s_1}} } \right)ds}  
	\\
	= \frac{{ - 1}}{{\Gamma \left( {2 + \alpha } \right)\Gamma \left( { - \alpha } \right)}}\left( {\left. {{{\left( {s - t} \right)}^{1 + \alpha }}\left( {\int\limits_s^\infty  {{{\left( {{s_1} - s} \right)}^{ - 1 - \alpha }}u\left( {{s_1}} \right)d{s_1}} } \right)} \right|_{s = t}^\infty  } \right.
	\\
	\left. { - \left( {1 + \alpha } \right)\int\limits_t^\infty  {{{\left( {s - t} \right)}^\alpha }\left( {\int\limits_s^\infty  {{{\left( {{s_1} - s} \right)}^{ - 1 - \alpha }}u\left( {{s_1}} \right)d{s_1}} } \right)ds} } \right) 
	\\
	= \frac{1}{{\Gamma \left( {1 + \alpha } \right)\Gamma \left( { - \alpha } \right)}}\int\limits_t^\infty  {{{\left( {s - t} \right)}^\alpha }\int\limits_s^\infty  {{{\left( {{s_1} - s} \right)}^{ - 1 - \alpha }}u\left( {{s_1}} \right)d{s_1}} ds}  
	\\
	= \frac{1}{{\Gamma \left( {1 + \alpha } \right)\Gamma \left( { - \alpha } \right)}}\int\limits_t^\infty  {\int\limits_t^{{s_1}} {{{\left( {s - t} \right)}^\alpha }{{\left( {{s_1} - s} \right)}^{ - 1 - \alpha }}ds\,} u\left( {{s_1}} \right)d{s_1}}  
	\\
	= \left| {\begin{array}{*{20}{c}}
			{s = \left( {{s_1} - t} \right)\tau  + t} \\ 
			{ds = \left( {{s_1} - t} \right)d\tau } 
	\end{array}} \right| 
	\\
	= \frac{1}{{\Gamma \left( {1 + \alpha } \right)\Gamma \left( { - \alpha } \right)}}\int\limits_t^\infty  {\int\limits_0^1 {{\tau ^\alpha }{{\left( {1 - \tau } \right)}^{ - 1 - \alpha }}d\tau \,} u\left( {{s_1}} \right)d{s_1}}  
	\\
	= \frac{{\Gamma \left( {1 + \alpha } \right)\Gamma \left( { - \alpha } \right)}}{{\Gamma \left( {1 + \alpha } \right)\Gamma \left( { - \alpha } \right)\Gamma \left( 1 \right)}}\int\limits_t^\infty  {u\left( {{s_1}} \right)d{s_1}}  = \int\limits_t^\infty  {u\left( s \right)ds}  = {}_tD_\infty ^{ - 1}u\left( t \right),
\end{multline*}
which is correct when 
\begin{equation}\label{VAS:S-B1.2}
	\mathop {\lim }\limits_{s \to \infty } \left( {{{\left( {s - t} \right)}^{1 + \alpha }}\left( {\int\limits_s^\infty  {{{\left( {{s_1} - s} \right)}^{ - 1 - \alpha }}u\left( {{s_1}} \right)d{s_1}} } \right)} \right) = 0.
\end{equation}

Equation \eqref{VAS:S-B1.1} reads then as 
\begin{equation}\label{VAS:S-B1.3}
	- {}_tD_\infty ^{ - 1}u(t) + A\,{}_tD_\infty ^{ - 2 - \alpha }u(t) = {}_tD_\infty ^{ - 2 - \alpha }f\left( t \right)
\end{equation}
and, when differentiating it by the rule 
\begin{multline*}
	\frac{d}{{dt}}{}_tD_\infty ^{ - 2 - \alpha }f\left( t \right) = \frac{1}{{\Gamma \left( {2 + \alpha } \right)}}\frac{d}{{dt}}\int\limits_t^\infty  {{{\left( {s - t} \right)}^{1 + \alpha }}f\left( s \right)ds}  \\
	=  - \frac{1}{{\Gamma \left( {1 + \alpha } \right)}}\int\limits_t^\infty  {{{\left( {s - t} \right)}^\alpha }f\left( s \right)ds}  =  - {}_tD_\infty ^{ - 1 - \alpha }f\left( t \right),
\end{multline*} 
we get the following for of \eqref{VAS:S-B1.3}, 
\begin{equation}\label{VAS:S-B1.4}
	u\left( t \right) - A{}_tD_\infty ^{ - \left( {1 + \alpha } \right)}u\left( t \right) =  - {}_tD_\infty ^{ - \left( {1 + \alpha } \right)}f\left( t \right).
\end{equation}

For the case $0 \leqslant \alpha  < 1$, we apply operator ${}_tD_\infty ^{ - 2 - \alpha }$ to the left and right hand-sides of the equation \eqref{VAS:S-B1.1}. For the first term, we will use the transformation
\begin{multline*}
	{}_tD_\infty ^{ - 2 - \alpha }{}_tD_\infty ^{1 + \alpha }u(t) \\
	=\frac{1}{{\Gamma \left( {2 + \alpha } \right)\Gamma \left( {1 - \alpha } \right)}}\int\limits_t^\infty  {{{\left( {s - t} \right)}^{1 + \alpha }}\frac{{{d^2}}}{{d{s^2}}}\left( {\int\limits_s^\infty  {{{\left( {{s_1} - s} \right)}^{ - \alpha }}u\left( {{s_1}} \right)d{s_1}} } \right)ds}  \\
	= \frac{1}{{\Gamma \left( {2 + \alpha } \right)\Gamma \left( {1 - \alpha } \right)}}\left( {\left. {{{\left( {s - t} \right)}^{1 + \alpha }}\frac{d}{{ds}}\left( {\int\limits_s^\infty  {{{\left( {{s_1} - s} \right)}^{ - \alpha }}u\left( {{s_1}} \right)d{s_1}} } \right)} \right|_{s = t}^\infty  - } \right.\\
	\left. { - \left( {1 + \alpha } \right)\int\limits_t^\infty  {{{\left( {s - t} \right)}^\alpha }\frac{d}{{ds}}\left( {\int\limits_s^\infty  {{{\left( {{s_1} - s} \right)}^{ - \alpha }}u\left( {{s_1}} \right)d{s_1}} } \right)ds} } \right) 
	\\
	= \frac{1}{{\Gamma \left( {2 + \alpha } \right)\Gamma \left( {1 - \alpha } \right)}}\left( {\left. { - \left( {1 + \alpha } \right)\left( {{{\left( {s - t} \right)}^\alpha }\left( {\int\limits_s^\infty  {{{\left( {{s_1} - s} \right)}^{ - \alpha }}u\left( {{s_1}} \right)d{s_1}} } \right)} \right)} \right|_{s = t}^\infty   } \right.
	\\
	\left. { + \alpha \left( {1 + \alpha } \right)\int\limits_t^\infty  {{{\left( {s - t} \right)}^{\alpha  - 1}}\int\limits_s^\infty  {{{\left( {{s_1} - s} \right)}^{ - \alpha }}u\left( {{s_1}} \right)d{s_1}} ds} } \right) 
	\\
	= \frac{1}{{\Gamma \left( \alpha  \right)\Gamma \left( {1 - \alpha } \right)}}\int\limits_t^\infty  {{{\left( {s - t} \right)}^{\alpha  - 1}}\int\limits_s^\infty  {{{\left( {{s_1} - s} \right)}^{ - \alpha }}u\left( {{s_1}} \right)d{s_1}} ds}  
	\\
	= \frac{1}{{\Gamma \left( \alpha  \right)\Gamma \left( {1 - \alpha } \right)}}\int\limits_t^\infty  {\int\limits_t^{{s_1}} {{{\left( {s - t} \right)}^{\alpha  - 1}}{{\left( {{s_1} - s} \right)}^{ - \alpha }}ds\,} u\left( {{s_1}} \right)d{s_1}}  \\
	= \left| {\begin{array}{*{20}{c}}
			{s = \left( {{s_1} - t} \right)\tau  + t} \\ 
			{ds = \left( {{s_1} - t} \right)d\tau } 
	\end{array}} \right| 
\end{multline*}
\begin{multline*}
	= \frac{1}{{\Gamma \left( \alpha  \right)\Gamma \left( {1 - \alpha } \right)}}\int\limits_t^\infty  {\int\limits_0^1 {{\tau ^{\alpha  - 1}}{{\left( {1 - \tau } \right)}^{ - \alpha }}d\tau \,} u\left( {{s_1}} \right)d{s_1}}  
	\\
	= \frac{{\Gamma \left( \alpha  \right)\Gamma \left( {1 - \alpha } \right)}}{{\Gamma \left( \alpha  \right)\Gamma \left( {1 - \alpha } \right)\Gamma \left( 1 \right)}}\int\limits_t^\infty  {u\left( {{s_1}} \right)d{s_1}}  = \int\limits_t^\infty  {u\left( s \right)ds}  = {}_tD_\infty ^{ - 1}u\left( t \right),
\end{multline*}
which is correct when the following condition is true 
\begin{equation}\label{VAS:S-B1.5}
	\begin{gathered}
		\mathop {\lim }\limits_{s \to \infty } \left( {{{\left( {s - t} \right)}^{1 + \alpha }}\frac{d}{{ds}}\left( {\int\limits_s^\infty  {{{\left( {{s_1} - s} \right)}^{ - \alpha }}u\left( {{s_1}} \right)d{s_1}} } \right)} \right) = 0, \hfill \\
		\mathop {\lim }\limits_{s \to \infty } \left( {{{\left( {s - t} \right)}^\alpha }\int\limits_s^\infty  {{{\left( {{s_1} - s} \right)}^{ - \alpha }}u\left( {{s_1}} \right)d{s_1}} } \right) = 0 \hfill 
	\end{gathered} 
\end{equation}
and, we again arrive at \eqref{VAS:S-B1.3} but after differentiation, to \eqref{VAS:S-B1.4}.

Let us show that \eqref{VAS:S-B1.1} follows from \eqref{VAS:S-B1.4}. To do that, we apply operator $\left( { - {}_tD_\infty ^{\left( {1 + \alpha } \right)}} \right)$ to the left and right-hand sides. 
Considering $0 < 1 + \alpha < 2,$ we will use the following formula
${}_tD_\infty ^{\left( {1 + \alpha } \right)}{}_tD_\infty ^{ - \left( {1 + \alpha } \right)}p\left( t \right) = p\left( t \right)$, to get 
\[ - D_\infty ^{\left( {1 + \alpha } \right)}u\left( t \right) + AD_\infty ^{\left( {1 + \alpha } \right)}{}_tD_\infty ^{ - \left( {1 + \alpha } \right)}u\left( t \right) = D_\infty ^{\left( {1 + \alpha } \right)}{}_tD_\infty ^{ - \left( {1 + \alpha } \right)}f\left( t \right),\] 
from where \eqref{VAS:S-B1.1} follows. 

Therefore, we proved the following lemma.

\begin{lemma} \label{VAS:LmS-B 1}
	The problem \eqref{VAS:S-B1.1} is equivalent to \eqref{VAS:S-B1.4}  in the case $-1<\alpha<0$ when condition \eqref{VAS:S-B1.2} is fulfilled and, for the case $0\le\alpha<1$, when condition \eqref{VAS:S-B1.5} is satisfied.
\end{lemma}

Let us employ a fixed-point iteration to study existence of a solution of  the problem  \eqref{VAS:S-B1.4}. Then
\begin{multline}\label{VAS:S-B1.7}
	u\left( t \right) = A{}_tD_\infty ^{ - \left( {1 + \alpha } \right)}u\left( t \right) - {}_tD_\infty ^{ - \left( {1 + \alpha } \right)}f\left( t \right) \\
	= A\,{}_tD_\infty ^{ - \left( {1 + \alpha } \right)}\left( {A{}_tD_\infty ^{ - \left( {1 + \alpha } \right)}u\left( t \right) - {}_tD_\infty ^{ - \left( {1 + \alpha } \right)}f\left( t \right)} \right) - {}_tD_\infty ^{ - \left( {1 + \alpha } \right)}f\left( t \right) \\
	= {A^2}\,{}_tD_\infty ^{ - 2\left( {1 + \alpha } \right)}u\left( t \right) - A\,{}_tD_\infty ^{ - 2\left( {1 + \alpha } \right)}f\left( t \right) - {}_tD_\infty ^{ - \left( {1 + \alpha } \right)}f\left( t \right) =  \ldots  \\
	= {A^n}\,{}_tD_\infty ^{ - n\left( {1 + \alpha } \right)}u\left( t \right) - \sum\limits_{k = 1}^n {{A^{k - 1}}\,{}_tD_\infty ^{ - k\left( {1 + \alpha } \right)}f\left( t \right)} .
\end{multline}

Introducing the notation $\nu  = \alpha  + 1$ with rational $\nu $, i.e., 
\begin{equation}\label{VAS:S-B1.8}
	\nu  = \frac{p}{q},\quad p,q \in \mathbb{N}
\end{equation}
taking $n = q$ in \eqref{VAS:S-B1.7} deduce 
\begin{multline*}
	u\left( t \right) = {A^q}\,{}_tD_\infty ^{ - p}u\left( t \right) - \sum\limits_{k = 1}^q {{A^{k - 1}}\,{}_tD_\infty ^{ - k\nu }f\left( t \right)}  =  \hfill \\
	= \frac{{{A^q}}}{{\Gamma \left( p \right)}}\,\int\limits_t^\infty  {{{\left( {s - t} \right)}^{p - 1}}u\left( s \right)ds}  - \sum\limits_{k = 1}^q {{A^{k - 1}}\,{}_tD_\infty ^{ - k\nu }f\left( t \right)}
\end{multline*}
whose differentiation $p$ times with respect to $t$ gives 
\begin{equation}\label{VAS:S-B1.9}
	\frac{{{d^p}u\left( t \right)}}{{d\,{t^p}}} = {\left( { - 1} \right)^p}{A^q}\,u\left( t \right) - \frac{{{d^p}}}{{d\,{t^p}}}\sum\limits_{k = 1}^q {{A^{k - 1}}\,{}_tD_\infty ^{ - k\nu }f\left( t \right)} .
\end{equation}

We consider the following equation
\begin{equation}\label{VAS:S-B1.10}
	\frac{{{d^p}y\left( t \right)}}{{d\,{t^p}}} = {\left( { - 1} \right)^p}{a^q}\,y\left( t \right) - F\left( t \right),
\end{equation}
where $a  $ is a constant, to study \eqref{VAS:S-B1.9}.

It is a linear equation of the $p$-th order whose solution can be presented as the sum of its `homogeneous and inhomogeneous' solutions. 
To determine both summands, the characteristic equation is decisive
\begin{equation}\label{VAS:S-B1.11}
	{\lambda ^p} = {\left( { - 1} \right)^p}{a^q}.
\end{equation}

In the case of the Cauchy problem for the homogeneous equation \eqref{VAS:S-B1.10}, as it was previously shown (in the previous paragraph), the solution can be written using an exponent with choosing the solution of the characteristic equation that provides decrease at infinity.
The solution of the inhomogeneous equation \eqref{VAS:S-B1.10} with homogeneous initial conditions can be written using the convolution integral
\begin{equation}\label{VAS:S-B1.12}
	y\left( t \right) =  - \int\limits_0^t {K\left( {t - s} \right)F\left( s \right)ds} ,
\end{equation}
where the kernel $K\left( t \right)$ satisfies the following Cauchy problem
\begin{equation}\label{VAS:S-B1.13}
	\begin{cases}
		\displaystyle {\frac{{{d^p}K\left( t \right)}}{{d\,{t^p}}} = {{\left( { - 1} \right)}^p}{a^q}\,K\left( t \right),} \\[2mm] 
		{K\left( 0 \right) = K'\left( 0 \right) =  \ldots  = {K^{(p - 2)}}\left( 0 \right) = 0,} \\[2mm] 
		{K^{(p - 1)}\left( 0 \right) = 1.} 
	\end{cases}  
\end{equation} 

Because the equation in \eqref{VAS:S-B1.13} is homogeneous, its solution can be found using \eqref{VAS:S-B1.11} whose roots are
\[{\lambda _k} = {a^{\frac{q}{p}}}\left( {\cos \left( {\frac{{\left( {p + 2k} \right)\pi }}{p}} \right) + i\,\sin \left( {\frac{{\left( {p + 2k} \right)\pi }}{p}} \right)} \right),\quad k = 0,1, \ldots ,p - 1;\] 
then $K\left( t \right)$ admits the form 
\begin{equation}\label{VAS:S-B1.14}
	K\left( t \right) = \sum\limits_{j = 0}^{p - 1} {{c_j}{\operatorname{e} ^{{\lambda _j}t}}},
\end{equation}
where constants ${c_j},\;j = \overline {0,p - 1} $ should come from the linear system (LS)
\begin{equation}\label{VAS:S-B1.15}
	\sum\limits_{j = 0}^{p - 1} {{\lambda _j}^k{c_j}}  = 0,\ k = \overline {0,p - 2};\ \ \ 
	\sum\limits_{j = 0}^{p - 1} {{\lambda _j}^{p - 1}{c_j}}  = 1,
\end{equation}
where 
with the Vandermonde matrix. 
Because the right-hand side has only one nonzero scalar element, 
the LS solution can be in the simplest way by constructed by the Kramer method,
\[
{c_k} = \frac{{\det \left( {{A_k}} \right)}}{{\det \left( A \right)}},\quad k = \overline {0,p - 1} ,
\] 
where $A$ is the LS matrix \eqref{VAS:S-B1.15},  ${A_k}$ is the  matrix $A$ with the $k$-th column replaced by the right-hand side vector. 
Expanding ${A_k}$  over the $k$-th column and accounting for that algebraic complement is also the Vandermonde determinant of order $\left( {p - 1} \right)$, the following solution is derived 
\begin{equation}\label{VAS:S-B1.16}
	{c_k} = \frac{{{{\left( { - 1} \right)}^{k + p}} \prod\limits_{ 0 \leqslant i < j \leqslant p - 1,  i,j \ne k }  {\left( {{\lambda _j} - {\lambda _i}} \right)} }}{{\prod\limits_{0 \leqslant i < j \leqslant p - 1} {\left( {{\lambda _j} - {\lambda _i}} \right)} }},\quad k = \overline {0,p - 1} .
\end{equation}
Since all ${\lambda _i}$ are different, it follows from \eqref{VAS:S-B1.16} that ${c_k},$ $k = \overline {0,p - 1} $ are not equal to 0. 
Therefore, it is important to study the sign of the real part of all ${\lambda _i}.$

\medskip 

Let us examine ${\lambda _i}.$ To do this, consider the possible cases\\
1)	$p = 1,$ then
$\lambda  =  - {a^q} < 0.$ \\
2)	$p = 2n,\;n \in {\rm N},$ then
\begin{multline*}
	{\lambda _k} = {a^{\frac{q}{{2n}}}}\left( {\cos \left( {\frac{{\left( {2n + 2k} \right)\pi }}{{2n}}} \right) + i\,\sin \left( {\frac{{\left( {2n + 2k} \right)\pi }}{{2n}}} \right)} \right) \\
	= {a^{\frac{q}{{2n}}}}\left( {\cos \left( {\frac{{\left( {n + k} \right)\pi }}{n}} \right) + i\,\sin \left( {\frac{{\left( {n + k} \right)\pi }}{n}} \right)} \right),\quad k = 0,1, \ldots ,2n - 1,
\end{multline*}
from where we get ${\lambda _n} = {a^{\frac{q}{{2n}}}}\left( {\cos \left( {2\pi } \right) + i\,\sin \left( {2\pi } \right)} \right) = {a^{\frac{q}{{2n}}}} > 0.$  \\
3)	$p = 2n + 1,\;n \in {\rm N},$ then
\begin{multline*}{\lambda _k} = {a^{\frac{q}{{2n + 1}}}}\left( {\cos \left( {\frac{{\left( {2n + 1 + 2k} \right)\pi }}{{2n + 1}}} \right) + i\,\sin \left( {\frac{{\left( {2n + 1 + 2k} \right)\pi }}{{2n + 1}}} \right)} \right),\\
	k = 0,1, \ldots ,2n.
\end{multline*}
As a result, 
\begin{multline*}
	{\lambda _n} = {a^{\frac{q}{{2n + 1}}}}\left( {\cos \left( {\frac{{\left( {4n + 1} \right)\pi }}{{2n + 1}}} \right) + i\,\sin \left( {\frac{{\left( {4n + 1} \right)\pi }}{{2n + 1}}} \right)} \right) \\
	= {a^{\frac{q}{{2n + 1}}}}\left( {\cos \left( {\left( {2 - \frac{1}{{2n + 1}}} \right)\pi } \right) + i\,\sin \left( {\left( {2 - \frac{1}{{2n + 1}}} \right)\pi } \right)} \right) \\
	= {a^{\frac{q}{{2n + 1}}}}\left( {\cos \left( {\frac{\pi }{{2n + 1}}} \right) - i\,\sin \left( {\frac{\pi }{{2n + 1}}} \right)} \right).
\end{multline*}
Because $n \in {\rm N},$ 
$\operatorname{Re} \left( {{\lambda _n}} \right) > 0$, thus, among solutions of the characteristic equation \eqref{VAS:S-B1.1} for $p \geqslant 2$ there exists a solution whose real part is positive. 
This means that in the case where there is an unbounded operator in \eqref{VAS:S-B1.9}, the solution representation \eqref{VAS:S-B1.12} has components with an operator exponent of the form
$\exp \left( {{A^{\frac{q}{p}}}a\,t} \right),$
where $a = const > 0.$ This means that function $F\left( t \right)$ is analytic one of operator ${A^{\frac{q}{p}}}.$  
Hence, the following theorem is true. 

\begin{theorem}\label{VAS:Th1-S-B}
	Let the conditions of Lemma \ref{VAS:LmS-B 1} be fulfilled. Then \eqref{VAS:S-B1.1} and its equivalent form \eqref{VAS:S-B1.4} have the only solution for the case 
	$
	\alpha={1}/{q}-1,
	$
	$q\in  \mathbb{N}.$
\end{theorem}


\subsection{Representation of the solution of inhomogeneous differential equations with fractional derivatives and unbounded operator coefficients in the Banach space}

According to Theorem \ref{VAS:Th1-S-B},
$$\alpha+1=\frac{1}{q},\quad q\in  \mathbb{N}$$
and, therefore, the problem \eqref{VAS:S-B1.4}, as it was shown above, reduces to the form \eqref{VAS:S-B1.9} with $p=1$, that is, 
\begin{equation}\label{VAS:S-B1.17}
	\frac{{{d}u\left( t \right)}}{{d\,{t}}} = -{A^q}\,u\left( t \right) - F(t),\ \ F(t)=\frac{{{d}}}{{d\,{t}}}\sum\limits_{k = 1}^q {{A^{k - 1}}\,{}_tD_\infty ^{ - k\nu }f\left( t \right)} 
\end{equation}
and, employing results of \cite{VAS:GMV-mon}, the solution representation for \eqref{VAS:S-B1.17} takes the form 
\begin{equation}\label{VAS:S-B1.18}
	u(t)={\rm e}^{-A^qt}u(0)- \int\limits_0^t {\rm e}^{-A^q(t-s)}F(s)ds=u_h(t)+u_{ih}(t).
\end{equation}

\subsection{Numerical method for the inhomogeneous equation}

For the first term in \eqref{VAS:S-B1.18}, we will use the approximate method from subsection \ref{VAS:metod-odn}.
To do that {$u_{ih}(t)$}, we apply the algorithm proposed in \cite{VAS:GMV-mon}, i.e., 
\begin{multline}\label{VAS:ih1-nab}
	u_{ih}(t)\approx u_{1,N}(t)=\frac{h}{2 \pi
		i}\sum_{k=-N}^{N}z^{\prime}(kh)[(z(kh)I-A)^{-1}-\frac{1}{z(kh)}I]
	\\
	\times h\sum_{p=-N}^{N}\mu_{k,p}(t)F(\omega_p(t)),
\end{multline}
where
$$
\mu_{k,p}(t)=\frac{t}{2\cosh^2{(ph)}}\text{exp}\{-\frac{t}{2} z(kh)[1-\tanh{(ph)}]\}, $$
$$
\omega_p(t)=\frac{t}{2}[1+\tanh{(ph)}]; \ \ h=\frac{1}{\sqrt{N}},
$$
$$
z(\xi)=a_I \cosh{\xi}-i b_I \sinh{\xi};\ \ \  z'(\xi)=a_I \sinh{\xi}-i
b_I \cosh{\xi}.
$$

The following theorem specifies convergence and error of the algorithm.

\begin{theorem}\label{VAS:inh1-conv}
	Let the conditions of Theorem \ref{VAS:th1} be fulfilled, $F(t) \in L((0;T),X)$ and $F(t) \in D(A ^{\sigma}),$ $\sigma>0$ $ \forall t \in [0,\infty]$ has an analytical continuation into the sector $$\Sigma_f=\{ d {\rm e}^{i \theta_1 }: \; d \in [0,\infty], \;|\theta_1|<\varphi\},$$ where the estimate
	\begin{equation}\label{VAS:analie24000}
		\|A^{\sigma}F(w)\|\le c_\alpha {\rm e}^{-\delta_\sigma |\text{Re}\; w|}, \; w \in \Sigma_f, \; \delta_\sigma \in (0,\sqrt{2}\rho]
	\end{equation}
	holds true. Then (\ref{VAS:ih1-nab}) converges to $u_{ih}(t)$ and the error estimate 
	\begin{equation}\label{VAS:ie250}
		\|{\cal E}_N(t)\|=\|u_{ih}(t)-u_{1,N}(t)\|\le c {\rm e}^{-c_1\sqrt{N}}
	\end{equation}
	uniformly over $t\ge 0$ under condition that $h={\cal O}(1/\sqrt{N})$. The constants $c$, $c_1$ are positive, depending on $\sigma,$ $\varphi,$ $\rho$ and independent of $N$, $t$.
\end{theorem}
Note that \eqref{VAS:ih1-nab} admits parallel computations. 
For each term, resolvents $$\left[(z(kh)I-A)^{-1}-\frac{1}{z(kh)}I \right],$$ can be independently computed. In addition, Theorem~\ref{VAS:inh1-conv} states that the method provides an exponential rate of the convergence.

%
%

\end{document}